\documentstyle[12pt]{article}

  \newcommand{\const}{\rm const}
  
  \newcommand{\supp}{\rm supp}

\begin{document}

   \begin{center}

 {\bf Classical non - linear operators} \par

 \vspace{4mm}

{\bf  in Grand Lebesgue Spaces.} \par

\vspace{5mm}

{\bf M.R.Formica,  E.Ostrovsky and L.Sirota}.\par

\end{center}

 Universit\`{a} degli Studi di Napoli Parthenope, via Generale Parisi 13, Palazzo Pacanowsky, 80132,
Napoli, Italy. \\

e-mail: mara.formica@uniparthenope.it \\

\vspace{4mm}

Department of Mathematics and Statistics, Bar-Ilan University, \\
59200, Ramat Gan, Israel. \\

e-mail: eugostrovsky@list.ru\\

\vspace{4mm}

\ Department of Mathematics and Statistics, Bar-Ilan University, \\
59200, Ramat Gan, Israel. \\

e-mail: sirota3@bezeqint.net \\

\vspace{5mm}

\begin{center}

{\bf Abstract.}

\end{center}

 \hspace{3mm} We study in this short report the boundedness of classical non - linear operators:
 Nemytskii, Urysohn, Hammerstein acting from  one Grand Lebesgue Space to another one, and deduce
 some its {\it upper} norm estimates. \par
  \ We bring also some examples to illustrate the exactness of our estimates. \par

\vspace{5mm}

  \hspace{3mm} {\it Key words and phrases.} Measure and measurable space,  Nemytskii, Urysohn, Hammerstein
and other operators, boundedness,  ordinary and operators norm, estimating factorization, power operator, degree,
ordinary and mixing (anisotropic)
Lebesgue - Riesz  and Grand Lebesgue Spaces (GLS), estimate, H\"older's inequality,  examples, measurable functions,
generating function. \par

\vspace{5mm}

\section{Statement of problem. Notations. Definitions.}

\vspace{5mm}

 \hspace{3mm} Let  $ \ (X = \{x\}, \ {\cal M}, \ \mu); \ (Y = \{y\}, {\cal N}, \nu);  \ $  be
 two measurable spaces equipped with two non - trivial measures  correspondingly  $ \ (\mu, \nu). \ $ \ Let also
 $ \ U = U[g](x), \ $  where

\begin{equation} \label{general oper}
f(x): = K[g(\cdot)](x) =  K[g(\cdot), x, \cdot],
\end{equation}

$$
 f: X \to R, \ g: Y \to R
$$
 be an {\it operator,} not necessary to be linear.\par

\vspace{3mm}

\begin{center}

   {\it It is presumed that all the considered functions are measurable.} \par

\end{center}

\vspace{5mm}

 \hspace{3mm} {\bf We intent in this short report to estimate the norms of these operators as an operator acting between two
 Grand Lebesgue Spaces (GLS). } \par

 \vspace{5mm}

  \ We will consider the following three  important types of nonlinear operators: Nemytskii, Urysohn and Hammerstein.\par

\vspace{5mm}

 \begin{center}

 {\sc Definitions of operators.}

 \end{center}

\vspace{4mm}

 \ {\bf Definition 1.} Nemytskii operator, see   \cite{Besov 2},  \cite{Krasnoselskii 2}, \cite{Nemytsky},  \cite{Urysohn}:

\begin{equation} \label{Nemytskii}
N[g](x) \stackrel{def}{=} n(x,g(x)).
\end{equation}

\vspace{4mm}

 \ {\bf Definition 2.} Urysohn's operator, see   \cite{Kantorowitch Akilov}, pp. 423 - 429, \cite{Krasnoselskii 1}, \cite{Urysohn}:

\begin{equation} \label{Urysohn}
U[g](x) \stackrel{def}{=}   \int_Y  u(x,y,g(y))  \ \nu(dy).
\end{equation}

\vspace{4mm}

 \ {\bf Definition 3.} Hammerstein's operator, \  \cite{Dolph},  \cite{Hammerstein}, \cite{Howard},  \cite{Krasnoselskii 2}:

\begin{equation} \label{Hammerstein}
H[g](x) \stackrel{def}{=} \int_Y h(x,y) \ w(y, g(y)) \ \nu(dy).
\end{equation}

\vspace{4mm}

 \ Here all the introduced functions $ \ n(x,y), \ u(x,y,z), \ h(x,y), \ w(x,y) \ $ are measurable and numerical valued. \par
 
\vspace{4mm}

 \ These operators appear in particular in the theory of non - linear  Integral and  Partial Differential Equations (PDE),
see e.g. \cite{Chen 1} - \cite{Chen 4}, \cite{Renardy},  \cite{Runst}, and reference therein. \par

\vspace{5mm}

 \ {\sc  A brief excursus in the theory of Grand Lebesgue Spaces (GLS). } \

\vspace{4mm}

\ Let $ \ (V = \{v\}, {\cal D}, \gamma)  \ $ be again measurable space equipped with non - trivial sigma finite measure $ \ \gamma. \ $ Let also
$ \ (a,b) = \const, \ 1 \le a < b \le \infty  \ $ and let $ \ \psi = \psi(p), \ a < p < b \ $ be certain strictly positive numerical valued function;
which is named as ordinary {\it generating function.}  By definition, the Grand Lebesgue Space (GLS)  $ \ G\psi[V,\psi; a,b] = G\psi \ $ consists on
all the numerical valued measurable functions $ \ f: V \to R  \ $ having the following finite norm

\begin{equation} \label{Gpsi norm}
||f||G\psi = ||f||G\psi[V,\psi; a,b] \stackrel{def}{=} \sup_{p \in (a,b)} \ \left\{ \ \frac{||f||_{p,V}}{\psi(p)} \ \right\}.
\end{equation}

 \ Hereafter $ \ ||f||_{p,V} \ $ denotes the usually Lebesgue - Riesz norm:

\begin{equation} \label{Lebesgue Riesz}
 ||f||_{p,V} \stackrel{def}{=} \left[ \ \int_V |f(v)|^p \ \gamma(dv) \ \right]^{1/p}, \ 1 \le p < \infty.
\end{equation}

\vspace{4mm}

 \ Notation: $ \ (a,b) \stackrel{def}{=} \ \supp \ \psi. \ $ \par
 \ Define formally $ \ \psi(p) = \infty,  \ $ when $ \ p \ \notin \ \supp \ \psi. \ $ \par

 \vspace{3mm}

  \ {\it Definition of the natural generating function.} \ Let $ \ f = f(v)  \ $  be measurable numerical valued function such that

$$
\exists (a,b), \ 1 \le a < b \le \infty \ \Rightarrow \forall p \in (a,b) \ ||f||_{p,V} < \infty.
$$

 \ The natural generating function $ \ \psi[f](p) \ $ is defined as follows:

$$
\psi[f](p) \stackrel{def}{=} ||f||_{p,V}, \ p \in (a,b).
$$

  \ Evidently, $ \ ||f||G\psi[f] = 1. \ $ \par

\vspace{4mm}

 \ The general theory of these spaces is represented in many works, see e.g. \cite{Ahmed Fiorenza Formica at all}, \cite{anatriellofiojmaa2015},
\cite{anatrielloformicaricmat2016},  \cite{Buldygin-Mushtary-Ostrovsky-Pushalsky}, \cite{caponeformicagiovanonlanal2013},  \cite{Ermakov etc. 1986},
\cite{Fiorenza2000}, \cite{fiokarazanalanwen2004}, \cite{fioguptajainstudiamath2008}, \cite{Fiorenza-Formica-Gogatishvili-DEA2018},
\cite{fioforgogakoparakoNAtoappear}, \cite{fioformicarakodie2017},\cite{Kozachenko}, \cite{Ostr Sir CLT  mixed},\cite{Ostrovsky 0},
\cite{Ostrovsky 1}, \cite{Ostrovsky HAIT} etc. In particular, these spaces are complete, Banach functional and rearrangement invariant.\par

\vspace{4mm}

 \ It is known, see e.g. \cite{Ostrovsky 2}, that in the case when $ \ a = 1, b = \infty \ $ the space $ \ G\psi_{1,\infty} \ $ coincides
with appropriate exponential Orlicz space. The belonging of some function $ \ h = h(v) \ $ is closely related with its tail behavior

$$
T_h(t) := \gamma \{ \ v: \ |h(v)| \ge t \ \},  \ t \to \infty.
$$

\vspace{5mm}

 \ \section{Main  result: Nemytskii operator.}

\vspace{5mm}

 \hspace{3mm} We consider here the Nemytskii's operator  (\ref{Nemytskii}). Let us impose the  following  condition
of factorization estimating on the leading function  $ \ n = n(x,y): \ $

 \begin{equation} \label{nem cond 1}
\exists \beta = \const \in [1,\infty), \ \exists \phi = \phi(x), \  x \in X, \ \Rightarrow |n(x,y)| \le \phi(x) \cdot |y|^{\beta},
 \end{equation}

 \ $ \ \phi: X \to R. \ $  Assume further that in  (\ref{Nemytskii})  [and in  (\ref{nem cond 1})]

\begin{equation} \label{cond factoriz est}
 g(\cdot)  \in G\psi, \  \phi(\cdot) \in G\nu,
\end{equation}
for suitable generating functions $ \ \psi, \ \nu; \ $  so that

$$
||g||_{p,X} \le C_1 \ \psi(p), \ ||\phi||_{q,X} \le C_2 \nu(q),
$$
where $ \ C_1 = ||g||G\psi, \ C_2 = ||\phi||G\nu. \ $  \par

 \ Of course, one can choose as such a functions $ \ \psi, \ \nu \ $ the natural ones: $ \ \psi(p) = ||g||_{p,X}, \ \nu(q) = ||\phi||_{q,X},\ $
if they there exist. Then $ \  C_1 = C_2 = 1. \ $ \par

\  Following,  if we introduce the so - called {\it power} operator of the degree $ \ \beta: \     P_{\beta} [g] (x) = |g(x)|^{\beta}, $ where as before
$ \ \beta = \const \ge 1, \ $ then

$$
||P_{\beta}[g]||_{p,X} = ||g^{\beta}||_{p,X} \le C_1^{\beta} \cdot \psi^{\beta}(\beta p).
$$

\vspace{3mm}

 \ We have:

$$
|f(x)| \le |\phi(x)| \cdot |g(x)|^{\beta}.
$$

 \ One can apply H\"older's inequality:

\begin{equation} \label{Holder}
||f||_{r,X} \le ||\phi||_{q,X} \cdot ||g^{\beta}||_{p,X},
\end{equation}

where $ \ (p,q,r) \ $ are arbitrary numbers such that $ \ p,q,r \ge 1 \ $ and moreover

$$
\frac{1}{r} =  \frac{1}{p} + \frac{1}{q},
$$
or equally

$$
q = \frac{pr}{p-r}, \hspace{3mm} 1 < r < p.
$$

\vspace{3mm}

 \ Introduce the following auxiliary functions

\begin{equation} \label{aux W}
W_a[\psi,\nu;\beta](p,r) =  W_a(p,r) \stackrel{def}{=}   \nu \left( \ \frac{pr}{p-r}  \ \right) \times \psi^{\beta}(\beta p),
\end{equation}

\begin{equation} \label{real W}
W[\psi, \nu; \ \beta](r) =   W(r) \stackrel{def}{=} \inf_{p > r} W_a(p,r).
\end{equation}

\vspace{3mm}

 \ It  follows from the inequality (\ref{Holder})

\begin{equation} \label{Hold estim}
||f||_{r,X} \le C_1^{\beta} \ C_2  \ W_a(p,r), \ p > r.
\end{equation}

\vspace{3mm}

 \  To summarize: \ we have from  (\ref{Hold estim}) after optimization over $ \ p; \ p \in(r, \ \infty): \ $ \par

\vspace{5mm}

 \ {\bf Theorem 2.1.} Suppose that for some non - trivial segment $ \  r \in (c,d): \ 1 \le c < d \le \infty \ $ the function
$ \ W[\psi,\nu; \beta](r) \ $ is finite. Then we propose in (\ref{Nemytskii}) under our notations

\begin{equation} \label{Main Nemytskii}
||f||GW[\psi,\nu; \beta]  \le  ||\phi||G\nu \times [ \ ||g||G\psi \ ]^{\beta}.
\end{equation}

\vspace{5mm}

  \ {\bf  Example 2.1.} Let us show by means of bringing examples the exactness of proposition (\ref{Main Nemytskii}), yet for arbitrary positive value
  $ \ \beta. \ $ Suppose $ \ \mu(X) = 1, \ \phi(x) = 1, \ $  so that the natural function for the $ \ \phi(\cdot) \ $ is equal to 1:

$$
\nu(q) = ||\phi||_{q,X} = 1, \ q \in [1,\infty).
$$

\vspace{3mm}

   \ Put also  $ \  \psi(p) = ||g||_p, \ 1 \le p < \infty, \ $ i.e. $ \ \psi(\cdot) \ $ is the natural
  generating function  for $ \ g(\cdot).\ $  Choose  for  definiteness  as a function $ \  g = g(x)  \  $  certain positive
  measurable bounded function having support with finite measure.\par
   \ Notice that this function $ \ p \to \psi(p), \ 1 \le p < \infty  \ $ is monotonically increasing. \par

 \vspace{3mm}

 \  On the other words,

  $$
  f(x) = g^{\beta}(x), \ x \in X.
  $$

 \ Let us investigate both the sides of the assertion of theorem 2.1 (\ref{Main Nemytskii}) in this case. We calculate the left hand side taking
into account the monotonicity and equality $ \ \nu(q) = 1 \ $ at the value $ \ p = r + 0 \ $

$$
W_a[\psi,\nu;\beta](r + 0,r) =  W_a(r + 0,r) = W(r) =   \psi^{\beta}(\beta r),
$$

 \ We have at the same time  for the right hand side of (\ref{Main Nemytskii}) choosing also $ \  p = r+0, \ $  so that $ \ q = + \infty:  \ $

$$
 ||g^{\beta}||_r  = \psi^{\beta}(\beta r), \ 1 \le r < \infty,
$$
which coincides with the right - hand side (in this example). \par

\vspace{5mm}

 \ \section{Main result. Urysohn's operator.}

\vspace{5mm}

 \hspace{3mm} We impose the following condition on the "kernel" $ \ u = u(x,y,z)  \ $ in  the definition of the Urysohn's operator
 (\ref{Urysohn}): \ $ \ \exists u_0 = u_0(x,y) \in R, \ x \in X, \ y \in Y,  \ $

\vspace{3mm}

\begin{equation} \label{cond Urysohn}
|u(x,y,z)| \le u_0(x,y) \ |z|^{\beta},  \ \exists  \beta = \const \ge 1.
\end{equation}

\vspace{3mm}

 \ Denote for at the same values $ \ (p,q,r) \ $ as before

\begin{equation} \label{kappa}
\kappa(q,r) \stackrel{def}{=} || \ || \ u_0 \ ||_{q,Y} \ ||_{r,X};
\end{equation}
here $ \ q = pr/(p-r), \ p > r > 1.  \ $ \par
  \ Assume as above that for some  non - trivial generating function $ \ \psi(\cdot) \ $

$$
  g(\cdot) \in G\psi.
$$

 \ We have by virtue of H\"older's inequality for these values of parameters $ \ (p,q,r) \ $

\begin{equation} \label{key Urysohn}
||U[g]||_{r,X} \le \left[ \ ||g||G\psi  \right]^{\beta}  \times [ \ \psi^{\beta}(\beta p) \ \kappa(pr/(p-r),r) \ ].
\end{equation}

\vspace{3mm}

 \ Introduce the following generation function

$$
\theta(r) \stackrel{def}{=} \inf_{p \in (r,\infty)} \left\{ \ \psi^{\beta}(\beta p) \times  \kappa  \left( \ \frac{pr}{p-r}, r \ \right) \ \right\}.
$$

 \ We conclude under our notations

 \vspace{5mm}

 \ {\bf  Theorem 3.1.}

\begin{equation} \label{Urysohn oper est}
||U[g]||G\theta  \le  \left[ \ ||g||G\psi \ \right]^{\beta}.
\end{equation}

 \vspace{5mm}

 \ {\bf Remark 3.1.} Unimprovability. The equality in the proposition of theorem 3.1 \ (\ref{Urysohn oper est}) \
 for all the positive values $ \ \beta \ $ is attained if for instance
  $ \ \nu(Y) = 1 = \mu(X), \ g(y) = C = \const \in (0,\infty), \ u_0(x,y) = 1, \ \ \kappa(q,r) = 1, \ \psi(p) = 1; \ $ then
 both the sides in (\ref{Urysohn oper est}) are equal to $ \ C^{\beta}. \ $ \par

\vspace{5mm}

 \ \section{Main result. Hammerstein's operator.}

\vspace{5mm}

 \hspace{3mm} The Hammerstein's operator, which we will rewrite as

\begin{equation} \label{Hammerstein rewrited}
H[g](x) \stackrel{def}{=} \int_Y h(x,y) \ n(y, g(y)) \ \nu(dy),
\end{equation}
 may be represented as a superposition of Nemytskii operator and Urysohn's one. Therefore, we
 retain all the notations and restrictions imposed above on the functions $ \ n(\cdot,\cdot). \ $ \par
 \ We have consequently by means of H\"older's inequality 

$$
 ||H[g]||_{r,X} \le || \ || \ h \ ||_{q,Y} \ ||_{r,X} \times ||n(\cdot, g(\cdot))||_{p,Y},
$$
 where as above $  \ 1/r = 1/q + 1/p, \ p,q,r > 1, \ p > r.  \ $  \par

 \vspace{3mm}

  \ The functional

\begin{equation} \label{aniz}
 h(\cdot,\cdot) \to || \ || \ h \ ||_{q,Y} \ ||_{r,X}   =  || \ || \ h \ || \ ||_{q,Y; r,X}
\end{equation}
 is named as ordinary {\it mixed}, or {\it anisotropic} Lebesgue - Riesz $ \ L_{q,r} \ $  norm of the function  $ \ h(\cdot,\cdot),  \ $
 see   \cite{Benedek Panzone}, \cite{Besov Ilin Nikolskii}, chapters 1,2. It is a multivariate generalization of the classical
 Lebesgue - Riesz norms.\par
  \ The correspondent multidimensional Grand Lebesgue Mixed (Anisotropic) norm  $ \ ||h(\cdot, \cdot)||G\tau   \ $ is defined very
  similar to the one - dimensional case

\begin{equation} \label{two dim Gtau}
 ||h(\cdot, \cdot)||G\tau  \stackrel{def}{=} \sup_{(q,r)}  \left\{ \ \frac{||h||_{q,r}}{\tau(q,r)}  \ \right\}.
\end{equation}

\vspace{3mm}

  \ Evidently, the space $ \ G \tau \  $ relative the norm $ \ ||\cdot||G\tau \ $ is complete bi - rearrangement invariant
Banach functional space. Herewith, for all the admissible values $ \ (q,r),  \ $  i.e. for which $ \ \tau(q,r) < \infty, \ $

\begin{equation} \label{two dim estim}
||h||_{q,r} \le ||h(\cdot, \cdot)||G\tau \cdot \tau(q,r).
\end{equation}

\vspace{4mm}

  \   Further,

$$
||n(\cdot, g(\cdot))||_{p,Y} \le ||\phi||_{s,Y}  \cdot ||g^{\beta}||_{t,Y}, \ 1/p = 1/r + 1/t, \ s = tp/(t-p), \ t > p.
$$
 \ Thus,

\begin{equation} \label{H eval}
||H[g]||_{r,X} \le \upsilon(r; p,t),
\end{equation}
where

\begin{equation} \label{upsilon}
 \upsilon(r; p,t):=  || \ || h \ ||_{q,Y} \ ||_{r,X} \times ||\phi||_{s,Y} \times ||g^{\beta}||_{t,Y}.
\end{equation}

\vspace{3mm}

 \ Introduce the following domain on the positive quadrant  $ \  D(r) = \{ \ p, t \   \}, \ $ dependent on the real number  (parameter)
 $ \ r; \ r  > 1: \ $ where

\begin{equation} \label{restric param}
\frac{1}{r} = \frac{1}{p} + \frac{1}{q}, \   \frac{1}{s} +  \frac{1}{t} =  \frac{1}{p}, \ p,q,r,s,t  > 1, \ p > r,\ t > p,
\end{equation}
and
$$
q = pr/(p - r), \ p > r; s = tp/(t - p),  \ t > p.
$$

\vspace{3mm}

 \ Denote also

\begin{equation} \label{Delta}
\Delta(r) := \inf_{(p,t) \in D(r)}  \upsilon(r; p,t).
\end{equation}

\vspace{3mm}

 \ Notice that if

\begin{equation} \label{three condit}
h(\cdot,\cdot) \in G\tau, \ \phi(\cdot) \in G\nu, \ g(\cdot) \in G\psi,
\end{equation}
then

\begin{equation} \label{ups estim}
\upsilon(r; p,t) \le ||h||G\tau \ ||\phi||G\nu \ [ \ ||g||G\psi \ ]^{\beta} \times \tau(q,r) \nu(s) \ \psi^{\beta}(\beta t).
\end{equation}

\vspace{4mm}

 \ We obtained actually the following proposition. \par

\vspace{4mm}

 \ {\bf Theorem 4.1.} We conclude in our notations $ \  H[g] \in G\Delta \ $   and moreover

\begin{equation} \label{main Hammerstein}
||H[g]||G\Delta \le 1.
\end{equation}

\vspace{5mm}

\vspace{0.5cm} \emph{Acknowledgement.} {\footnotesize The first
author has been partially supported by the Gruppo Nazionale per
l'Analisi Matematica, la Probabilit\`a e le loro Applicazioni
(GNAMPA) of the Istituto Nazionale di Alta Matematica (INdAM) and by
Universit\`a degli Studi di Napoli Parthenope through the project
\lq\lq sostegno alla Ricerca individuale\rq\rq .\par


\begin{thebibliography}{55}


\bibitem{Ahmed Fiorenza Formica at all}
{\bf I. Ahmed, A. Fiorenza, M.R. Formica, A. Gogatishvili, J.M. Rakotoson.}
{\it Some new results related to Lorentz G-Gamma spaces and interpolation.}
J. Math. Anal. Appl., 483, {\bf 2}, (2020).

\bibitem{Adams}
{\bf Adams R.A.}   {\it Anisotropic Sobolev Inequalities.} Casopic pro Pestovani Matematiky, (Prague), No. 3, 267 \ - \ 279.

\bibitem{anatriellofiojmaa2015}
{\bf G.~Anatriello} and {\bf A.~Fiorenza.} {\it Fully measurable
grand Lebesgue spaces}. J. Math. Anal. Appl. \textbf{422} (2015),
no.~2, 783--797.

\bibitem{anatrielloformicaricmat2016}
{\bf G.~Anatriello} and {\bf M.~R.~Formica.} {\it Weighted fully
measurable grand Lebesgue spaces and the maximal theorem}. Ric. Mat.
\textbf{65} (2016), no.~1, 221--233.

\bibitem{Arnold}
{\bf L.Arnold.} {\it  Stochastic Differential Equation; Theory and Applications.} Wiley, New York, 1974.

\bibitem{Benedek Panzone}
{\bf Benedek A. and Panzone R.} {\it The space} $ \ Lp \ $  {\it with mixed norm. Duke Math.}
J., 28, (1961), 301 - 324.

\bibitem{Besov Ilin Nikolskii}
{\bf Besov O.V., Ilin V.P., Nikolskii S.M.} {\it Integral representation of functions
and imbedding theorems. } Vol.1; Scripta Series in Math., V.H.Winston and Sons,
(1979), New York, Toronto, Ontario, London.


\bibitem{Besov 2}
{\bf  K.O.Besov.}  {\it On the Continuity of the Generalized Nemytskii Operator on Spaces of Differentiable Functions.}
 Mat. Zametki,2002, Volume 71, Issue 2, 168–181DOI: https://doi.org/10.4213/mzm337



\bibitem{Buldygin-Mushtary-Ostrovsky-Pushalsky}
{\bf V.~V.~Buldygin, D.~I.~Mushtary, E.~I.~Ostrovsky} and {\bf M.~I.~Pushalsky.} {\it New Trends in
Probability Theory and Statistics.} Mokslas (1992), V.1, 78--92;
Amsterdam, Utrecht, New York, Tokyo.


\bibitem{Chen 1}
{\bf Chen Wenxiong, Li Congming.} {\it Classification of solutions of some nonlinear elliptic
equations.} Duke Math. J., {\bf 63(3);}  \ 615 \ - \ 622, \ 1991.


\bibitem{Chen 2}
{\bf Chen Wenxiong, Li Congming.} {\it Classification of positive solutions for nonlinear differential
and integral systems with critical exponents.}  Acta Math. Sci. Ser. B (Engl.Ed.), {\bf 29(4):}, 949 \ - \ 960, 2009.


\bibitem{Chen 3}
{\bf Huang Genggeng. Li Congming.}  {\it A Liouville theorem for high order degenerate elliptic
equations.}  J. Differential Equations, {\bf 258 (4);} 1229 \ - \ 1251, \ 2015.

\bibitem{Chen 4}
{\bf Chen Dezhong, Ma Li.}  {\it  A Liouville type theorem for an integral system.}  Commun.
Pure Appl. Anal.,  {\bf 5(4); } \ 855 \ - \ 859, \ 2006.


\bibitem{Dunford Schwartz}
{\bf N.Dunford, B.Schwartz.} {\it Linear Operators. V.1, General Theory. } Academic Press (1958), New York,
London.

\bibitem{caponeformicagiovanonlanal2013}
{\bf C.~Capone, M.~R.~Formica} and {\bf R.~Giova.} {\it Grand
{L}ebesgue spaces with respect to measurable functions}. Nonlinear
Anal. \textbf{85} (2013), 125--131.


\bibitem{Dolph}
{\bf  C.L.Dolph.}  {\it Nonlinear Integral Equations of the Hammerstein Type.}
 Transactions of the American Mathematical Society,
Vol. 66, No. 2, (Jul., 1949), pp. 289-307


\bibitem{Dorogovtsev}
{\bf  A.A.Dorogovtsev.} {\it On applications of Gaussian random operator to random elements.}
Probab. Theory Appl., {\bf 30,} 1986, 812 \ - \ 814.


\bibitem{Ermakov etc. 1986}
{\bf S. V. Ermakov, and E. I. Ostrovsky.} {\it Continuity Conditions, Exponential Estimates, and the Central Limit Theorem for Random Fields.}
 Moscow, VINITY,  1986. (in Russian).

\bibitem{Fiorenza2000} {\bf A.~Fiorenza.} {\it Duality and reflexivity in grand Lebesgue
spaces.} Collect. Math. \textbf{51} (2000), no. 2, 131--148.

\bibitem{fiokarazanalanwen2004}
{\bf A.~Fiorenza} and {\bf G.~E.~Karadzhov.} {\it Grand and small
Lebesgue spaces and their analogs}, Z. Anal. Anwendungen \textbf{23}
(2004), no.~4, 657--681.

\bibitem{fioguptajainstudiamath2008}
{\bf A.~Fiorenza, B.~Gupta} and {\bf P.~Jain.} {\it The maximal
theorem for weighted grand Lebesgue spaces}. Studia Math.
\textbf{188} (2008), no.~2, 123--133.


\bibitem{Fiorenza-Formica-Gogatishvili-DEA2018}
{\bf A.~Fiorenza, M.~R.~Formica} and {\bf A.~Gogatishvili.} {\it On
grand and small Lebesgue and Sobolev spaces and some applications to
PDE's}. \emph{Differ. Equ. Appl.} \textbf{10} (2018), no.~1, 21--46.

\bibitem{fioforgogakoparakoNAtoappear}
{\bf A.~Fiorenza, M. R.~Formica, A.~Gogatishvili, T.~Kopaliani} and
{\bf J.~M. Rakotoson.} {\it Characterization of interpolation
between grand, small or classical Lebesgue spaces}. Preprint
arXiv:1709.05892, Nonlinear Anal., {to appear}.

\bibitem{fioformicarakodie2017}
{\bf A.~Fiorenza, M.~R.~Formica} and {\bf J.~M. Rakotoson.} {\it
Pointwise estimates for {$G\Gamma$}-functions and applications}.
Differential Integral Equations \textbf{30} (2017), no.~11-12,
809--824.


\bibitem{Hammerstein}
{\bf A.Hammerstein.}   {\it Nichtlineare Integralgleichungen nebst Anwendungen.} Acta Math. , {\bf 54,}
 (1930), pp. 117   \ - \ 176.


\bibitem{Howard}
{\bf Ralph Howard, Anton Schep.}  {\it  Norms of positive operators on} $ \ L^p \ $ {\it spaces. }

\bibitem{Kantorowitch Akilov}
{\bf L.V.Kantorowitch, G.P.Akilov.}  {\it Functional Analysis.} Moskow, "Science", 1984.


\bibitem{Krasnoselskii 1}
{\bf  M.A.Krasniselskii, Ya.B.Rutitskii.} {\it Convex functions aned Orlicz spaces.} GIFML, Moskow, 1958. (in Russian).

\bibitem{Krasnoselskii 2}
{\bf  M.A. Krasnosel'skii.}   {\it Topological methods in the theory of nonlinear integral equations.} , Pergamon,
(1964) (Translated from Russian).


\bibitem{Kozachenko}
{\bf Kozachenko Yu. V., Ostrovsky E.I.} (1985). {\it The Banach Spaces of random
Variables of subgaussian type.} Theory of Probab. and Math. Stat. (in
Russian). Kiev, KSU, 32, 43 - 57.

\bibitem{Ledoux Talagrand}
{\bf M.Ledoux, M.Talagrand.} {\it Probability in Banach spaces.} Springer Verlag, (1991).

\bibitem{Leoni}
{\bf D.Leoni.} {\it A first course on Sobolev spaces.}
Graduate Studies in Mathematics. Volume: 105; 2009; 607 pp;
MSC: Primary 46; Secondary 26.


\bibitem{Nemytsky}
{\bf Nemytsky V.V. }  {\it Theorems of existings and uniqueness  for nonlinear integral equations.}  Math. Sbornik, {\bf 4,}
  (1934),  (in Russian).

\bibitem{Okikiolu}
{\bf G.O.Okikiolu.}  {\it Aspects of the Theory of Bounded Integral Operators in Lp Spaces.} Academic Press,
(1971), London, New York.

\bibitem{Ostr Sir CLT  mixed}
{\bf E.Ostrovsky, L.Sirota.}   {\it Central Limit Theorem and exponential estimations in mixed (anisotropic)
Lebesgue Spaces.} \par
arXiv:1308.5606v1 [math.PR] 26 Aug 2013

\bibitem{Ostrovsky 0}
{\bf Ostrovsky E.I.}  (1999). {\it Exponential estimations for random Fields and its
Applications, (in Russian).}  Moscow - Obninsk, OINPE.

\bibitem{Ostrovsky 1}
{\bf E.Ostrovsky, L.Sirota and E.Rogover.} {\it Integral operators in bilateral Grand Lebesgue Spaces.} \par
arXiv:0912.2538v1 [math.FA] 13 Dec 2009

\bibitem{Ostrovsky HAIT}
{\bf E. Ostrovsky and L.Sirota.}  {\it Moment Banach spaces: theory and applications.}  HAIT Journal of Science
and Engeneering, C, Volume 4, Issues 1 - 2, pp. 233 - 262, (2007).


\bibitem{Ostrovsky 2}
{\bf Kozachenko Yu.V., Ostrovsky E., Sirota L.} {\it Relations between exponential tails, moments and
moment generating functions for random variables and vectors.} \\
arXiv:1701.01901v1 [math.FA] 8 Jan 2017


\bibitem{Ostrovsky}
{\bf E.Ostrovsky, L.Sirota} {\it Central Limit Theorem and exponential tail estimations in mixed
(anisotropic) Lebesgue spaces.} \\
arXiv:1308.5606v1 [math.PR] 26 Aug 2013


\bibitem{Renardy}
{\bf Renardy, Michael and Rogers, Robert C.}  (2004). {\it An introduction to partial differential equations.} Texts in
 Applied Mathematics 13 (Second ed.). New York: Springer-Verlag. p. 370. ISBN 0-387-00444-0. (Section 10.3.4)

\bibitem{Runst}
{\bf  RUNST, T. and SICKEL, W..} {\it Sobolev Spaces of Fractional Order, Nemytskii Operators and Nonlinear Partial
Differential Equations.}  De Gruyter, Berlin, 1996.

\bibitem{Urysohn}
{\bf P.S.Urysohn.}  {\it On the type of non - linear integral equations.}  Mathem. Sbornik, {\bf 36}, (1936), 213 \ - \  226,
(in Russian).


\end{thebibliography}
\end{document}